\documentclass[11pt,a4paper]{article}
\usepackage[utf8]{inputenc}
\usepackage[T1]{fontenc}
\usepackage{amsmath, amssymb, amsthm}
\usepackage{geometry}
\usepackage{graphicx}
\usepackage{cite}
\usepackage{float}
\usepackage{placeins}

\title{\textbf{Prime-Weighted Interference Patterns Inspired by the Euler Product}}
\author{Jouni J. Takalo}
\date{February 2026}

\newtheorem{definition}{Definition}
\newtheorem{proposition}{Proposition}
\newtheorem{remark}{Remark}

\begin{document}

\maketitle
\thispagestyle{empty} 

\begin{abstract}
We study a prime-weighted oscillatory model inspired by structural aspects of the Euler product of the Riemann zeta function. 
The model defines finite superpositions of prime-frequency modes and exhibits zero-like crossings produced by destructive interference. 
We analyze how the weight exponent $x$ controls amplitude growth, slope scaling, and stability of crossings. 
A heuristic asymptotic argument identifies $x=\tfrac12$ as a distinguished balance regime separating high-energy and over-damped behavior. 
The results concern the defined model itself.
\end{abstract}
\vspace{0.5em}

\clearpage
\setcounter{page}{1}
\newpage

\section{Introduction}

Prime-weighted trigonometric sums with logarithmic frequencies arise naturally from the Euler product representation of the Riemann zeta function (see e.g. \cite{Edwards_1974,Titchmarsh_1986}). When interpreted along vertical lines in the complex plane, such expressions generate oscillatory structures whose frequency spectrum is indexed by the prime numbers. Related oscillatory viewpoints appear in analytic number theory, in spectral interpretations of explicit formulae, and in connections to quantum-chaotic analogies \cite{Montgomery_1973,Odlyzko_1992 ,Berry_1986}.

In this work we isolate and study a finite prime-weighted cosine ensemble as a self-contained oscillatory system. Rather than investigating analytic properties of the Riemann zeta function $\zeta(s)$ itself, we analyze structural and scaling properties of the finite superpositions introduced in Definition~1 together with their phase-referenced variants (Definition~2). Our objective is to understand how interference mechanisms emerge in such ensembles and how the weight exponent $x$ governs amplitude growth, slope scaling, and the formation of zero-like crossings.

The central question may be formulated as follows:

\medskip
\noindent
\textit{How does the exponent $x$ control the balance between high-frequency richness and amplitude stability in finite prime-frequency superpositions?}
\medskip

We show that the model exhibits three robust regimes depending on whether $x < \tfrac12$, $x = \tfrac12$, or $x > \tfrac12$. A heuristic asymptotic analysis based on prime-density estimates suggests that $x = \tfrac12$ forms a distinguished balance exponent separating polynomial energy growth from bounded behavior. In this regime, destructive-interference minima sharpen with increasing prime cutoff while overall amplitude growth remains controlled.

Our approach is deliberately finite and structural. We do not attempt to approximate $\zeta(s)$, nor do we claim implications for the distribution of its zeros. Instead, the goal is to analyze intrinsic interference properties of a prime-frequency ensemble inspired by, but independent of, the Euler product framework.

The results combine heuristic asymptotic scaling arguments with numerical illustrations that visualize progressive prime superposition and regime transitions. Together they reveal a coherent scaling structure governed by the single parameter $x$.

\section{Model Definitions}

\begin{definition}[Raw prime cosine signal]
For a prime cutoff $P$ and exponent $x>0$, define
\[
S_{P,x}(t)=\sum_{p\le P} p^{-x}\cos(t\log p).
\]
\end{definition}

\begin{definition}[Phase-referenced signal]
For a prime cutoff $N$, exponent $x>0$, and phase reference $\theta(t)$, define
\[
W_{N,x}(t)
=
2\sum_{p\le N} p^{-x}\cos(\theta(t)-t\log p).
\]
\end{definition}

\begin{remark}
The signal $S_{P,x}(t)$ is used for quantitative scaling analysis and numerical experiments.
The phase-referenced signal $W_{N,x}(t)$ is used to illustrate phase-alignment mechanisms.
\end{remark}

\begin{definition}[Zero-like crossing]
A zero-like crossing is a point $t_0$ such that the signal vanishes with nonzero local slope.
\end{definition}
Intuitively, a zero-like crossing represents a true sign change of the signal. The nonvanishing slope condition ensures that the graph locally passes through zero rather than merely touching it. In the prime-weighted superposition, such crossings occur when destructive interference becomes strong enough to push the signal across the axis.

\section{Mechanism of Destructive Interference}

For each prime $p$, the component $\cos(t\log p)$ attains minima when
\[
t\log p \approx (2k-1)\pi.
\]
Different primes do not share identical minima, since their frequencies $\log p$ differ.
However, near-coincidence of many such minima within a short interval produces a pronounced negative excursion in the superposition.

\begin{proposition}[Heuristic zero-like crossing mechanism]
If for some $t_0$ a large subset of primes satisfies
\[
t_0\log p \approx (2k_p-1)\pi,
\]
then $S_{P,x}(t)$ develops a deep destructive-interference well near $t_0$, potentially producing a zero-like crossing.
\end{proposition}

\section{Amplitude and Balance Regimes}

Define the squared-amplitude budget
\[
B_P(x)=\sum_{p\le P} p^{-2x}.
\]

Using the prime number theorem heuristic
\[
\sum_{p\le P}f(p)\approx\int_2^P\frac{f(u)}{\log u}\,du,
\]
one obtains
\[
B_P(x)\approx\int_2^P\frac{u^{-2x}}{\log u}\,du.
\]

\begin{proposition}[Balance exponent]\label{prop:balance}
As $P\to\infty$,
\[
B_P(x)
\begin{cases}
\text{grows polynomially in }P, & x<\tfrac12,\\
\sim \log\log P, & x=\tfrac12,\\
\text{remains bounded}, & x>\tfrac12.
\end{cases}
\]
\end{proposition}

Thus $x=\tfrac12$ separates a high-energy regime from an over-damped regime and forms a distinguished balance exponent for the model.

\section{Slope Scaling}

Differentiating,
\[
S'_{P,x}(t)
=
-\sum_{p\le P} p^{-x}(\log p)\sin(t\log p).
\]

\begin{proposition}[RMS slope scaling (heuristic)]
Assuming quasi-random phase behavior,
\[
\mathbb{E}|S'_{P,x}(t)|^2
\approx
\frac12\sum_{p\le P}p^{-2x}(\log p)^2
\approx
\frac12\int_2^P u^{-2x}\log u\,du.
\]
Hence
\[
\big(\mathbb{E}|S'_{P,x}(t)|^2\big)^{1/2}
\asymp
\begin{cases}
P^{\frac{1-2x}{2}}\sqrt{\log P}, & x<\tfrac12,\\
\log P, & x=\tfrac12,\\
1, & x>\tfrac12.
\end{cases}
\]
\end{proposition}

This explains why crossings sharpen with increasing $P$ in the balance regime while remaining controlled.
\FloatBarrier

\section{Numerical Experiments}

\subsection{Phase alignment illustration}

\begin{figure}[H]
\centering
\includegraphics[width=0.8\textwidth]{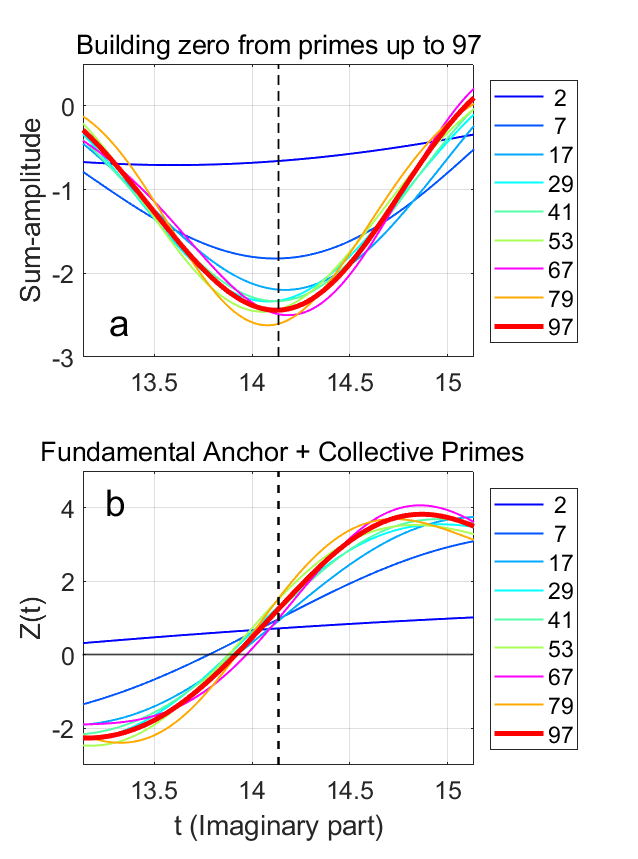}
\caption{
Progressive superposition of prime components in the phase-referenced signal $W_{N,1/2}(t)$.
(a) Partial sums obtained by including primes incrementally from $p=2$ to $p=97$.
Individual cosine components attain their minima at distinct phase locations determined by their frequencies $\log p$.
Although these minima do not coincide exactly, clusters of near-minima generate pronounced destructive-interference wells.
(b) The same process viewed under the phase reference $\theta(t)$, illustrating how increasing the prime cutoff sharpens the interference minimum and may produce a zero-like crossing.
}
\end{figure}

\subsection{Weight exponent comparison}

\begin{remark}
Figure 2 uses $S_{P,x}(t)$ with $t\in[140,160]$ sampled at 3000 points.
Two prime cutoffs are used: $P=100$ (25 primes) and $P=10^6$ (78{,}498 primes).
Red dashed lines indicate known zeta-zero ordinates for visual reference only.
\end{remark}

\begin{figure}[H]
\centering
\includegraphics[width=1.0\textwidth]{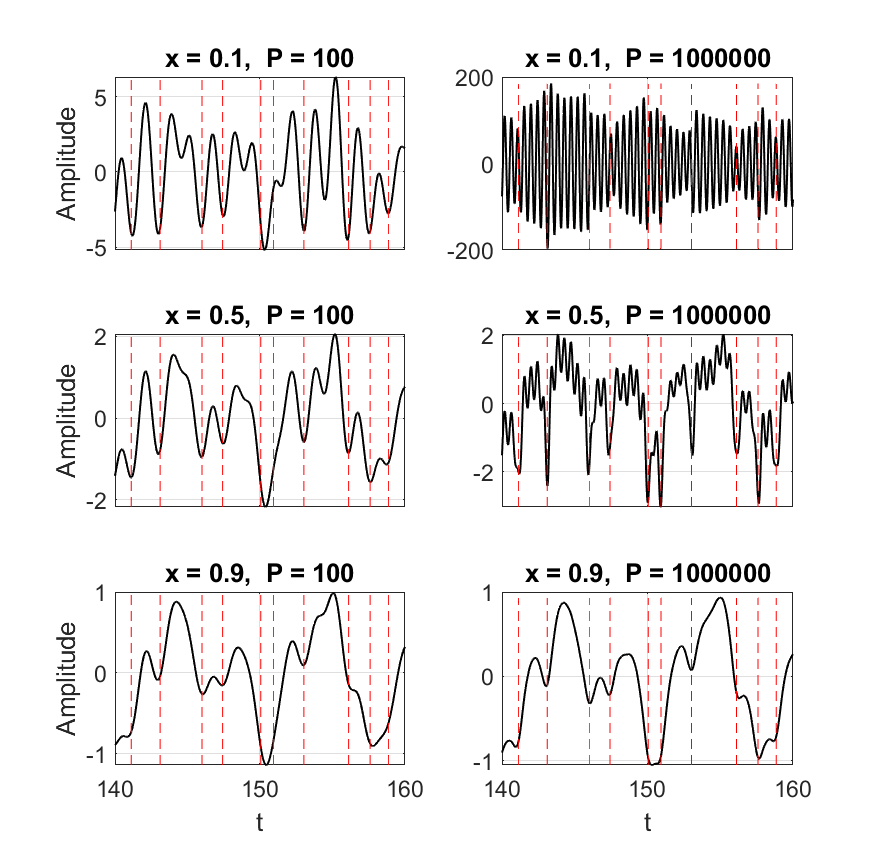}
\label{fig:stability}
\caption{
Weight-exponent dependence of the raw prime cosine signal 
$S_{P,x}(t)=\sum_{p\le P} p^{-x}\cos(t\log p)$
on $t\in[140,160]$.
Left: $P=100$ (25 primes).
Right: $P=10^6$ (78{,}498 primes).
Red dashed lines mark known zeta-zero ordinates (external reference).\\
\smallskip
\small
For $x<\tfrac12$, higher primes retain substantial weight and the large-ensemble signal becomes high-energy and rapidly oscillatory.\\
For $x=\tfrac12$, destructive-interference minima sharpen with increasing $P$ while overall amplitude growth remains controlled.\\
For $x>\tfrac12$, higher prime contributions are strongly attenuated, producing a smoother and effectively over-damped signal.}
\end{figure}

\section{Discussion}

The finite prime-weighted model exhibits three structurally distinct regimes governed by the exponent $x$. The transition at $x = \tfrac12$ is reflected simultaneously in amplitude scaling, RMS slope growth, and qualitative interference behavior.

For $x < \tfrac12$, higher primes retain substantial weight, leading to rapid polynomial growth of the squared-amplitude budget. The ensemble becomes increasingly high-energy and highly oscillatory as the prime cutoff increases. In this regime, interference minima occur but are unstable under extension of the prime range.

For $x > \tfrac12$, prime contributions decay sufficiently fast to yield bounded amplitude and controlled slopes. The resulting signal is comparatively smooth and effectively over-damped: higher primes contribute little structural refinement.

The balance case $x = \tfrac12$ occupies an intermediate position. Here the amplitude growth is logarithmic in scale, while slope growth remains sufficiently strong to sharpen interference minima as $P$ increases. This produces a regime in which destructive-interference wells become more sharply resolved without destabilizing global amplitude. The coexistence of structural richness and controlled energy distinguishes this exponent.

Several limitations should be emphasized. First, all asymptotic estimates are heuristic and rely on prime number theorem approximations together with quasi-random phase assumptions. Second, the analysis concerns finite ensembles and does not address analytic continuation or properties of $\zeta(s)$. The visual proximity between interference minima and known zeta-zero ordinates in numerical illustrations is presented solely for qualitative comparison and does not imply approximation results.

Nevertheless, the model suggests that prime-frequency ensembles possess an intrinsic scaling transition governed by $x$. From a broader perspective, the exponent $x = \tfrac12$ appears as the unique scaling at which frequency richness and amplitude control coexist in a nontrivial manner.

Possible directions for further study include probabilistic modeling of phase distributions, continuous prime-density analogues, spectral interpretations of the slope-scaling law, and investigation of stability properties under randomized prime perturbations. These questions lie beyond the scope of the present work but indicate that finite prime-weighted interference models may provide a useful framework for studying structured oscillatory systems indexed by arithmetic data.

\section{Conclusion}

We analyzed a prime-weighted interference model and identified $x=\tfrac12$ as a distinguished balance exponent separating high-energy and over-damped behavior.
Amplitude and slope scaling provide heuristic explanation for the observed regime transition.

\section*{Acknowledgments}

No personal or institutional support was involved.

\section*{Disclosure statement}

The author declares no conflicts of interest.

\section*{Data availability statement}

The numerical data supporting the findings of this study are available
from the author upon reasonable request.
\bibliographystyle{unsrt}
\bibliography{references_JT}

\end{document}